\RequirePackage{setspace}
\documentclass[aap,preprint]{imsart}
\usepackage{amsmath,amssymb,color}
\usepackage[T1]{fontenc}
\usepackage[francais]{babel}
\usepackage{paralist}
\usepackage{mathrsfs}
\usepackage{fancyhdr}
\usepackage{t1enc}
\usepackage[utf8]{inputenc}  
\usepackage[T1]{fontenc}    
\usepackage{hyperref}
\usepackage{cases}
\usepackage{setspace}
\usepackage[T1]{fontenc}      % Police contenant les caractères français
\usepackage[francais]{babel}  % Placez ici une liste de langues, la
\usepackage{paralist}
\usepackage{mathrsfs}
\usepackage[T1]{fontenc}
\usepackage{paralist}
\usepackage{mathrsfs}
\usepackage{amssymb}
\usepackage{lmodern}
\usepackage{amsmath}
\usepackage{setspace}
\usepackage{t1enc,amssymb,amsbsy}%ltimebuf}
\newcommand{\dproof}{\noindent {Proof.} \quad}
\newcommand{\fproof}{\hfill $\square$ \bigskip}
\startlocaldefs
\newtheorem{theorem}{Theorem}[part]
\newtheorem{definition}{Definition}[part]
\newtheorem{Proposition}{Proposition}[part]

\newtheorem{lemma}{Lemma}[part]
\newtheorem{corollary}{Corollary}[part]
\newtheorem{remark}{Remark}[part]
\newcommand{\cf}{\mathcal{F}}

\newcommand{\stops}{\mathcal{T}_{S}^p}
\newcommand{\stopo}{\mathcal{T}_0^p}
\newcommand{\stopsp}{\mathcal{T}_{S^+}^p}

\DeclareMathOperator*{\esssup}{ess\,sup} 
\DeclareMathOperator*{\essinf}{ess\,inf} 
\newcommand{\limsupn}{\limsup_{n\rightarrow\infty}}
\newcommand{\limn}{\lim_{n\rightarrow\infty}}

\endlocaldefs

\begin{document}

% "Title of the paper"
%\title{Optimal Stopping in General Predictable Framework}
%\runtitle{Optimal Stopping in General Predictable Framewor}

% indicate corresponding author with \corref{}
% \author{\fnms{John} \snm{Smith}\corref{}\ead[label=e1]{smith@foo.com}\thanksref{t1}}
% \thankstext{t1}{Thanks to somebody} 
% \address{line 1\\ line 2\\ printead{e1}}
% \affiliation{Some University}

%\author{\fnms{???} \snm{???}\ead[label=e1]{???}}
%\address{\printead{e1}}
%\and
%\author{\fnms{???} \snm{???}\ead[label=e2]{???}}
%\address{\printead{e2}}
%\affiliation{???}

%\runauthor{???}

%\begin{abstract}
%\end{abstract}

%\begin{keyword}[class=MSC]
%\kwd[Primary ]{}
%\kwd{}
%\kwd[; secondary ]{}
%\end{keyword}

%\begin{keyword}
%\kwd{}
%\kwd{}
%\end{keyword}
\begin{frontmatter}
	\title{Optimal Stopping in General Predictable Framework}
	\runtitle{Optimal Stopping in General Predictable Framework}
	%\thankstext{T1}{Footnote to the title with the ``thankstext'' command.}
	
	\begin{aug}
		\author{\fnms{Siham} \snm{Bouhadou}\thanksref{m1}\ead[label=e1]{sihambouhadou@gmail.com}}
		\and
		\author{\fnms{Youssef} \snm{Ouknine}\thanksref{m1,m2}\ead[label=e2]{ouknine@uca.ac.ma} \ead[label=e3]{youssef.ouknine@um6p.ma}}

\runauthor{ Bouhadou and Ouknine}

		\affiliation{Cadi Ayyad University\thanksmark{m1} and Mohammed VI Polytechnic University\thanksmark{m2}}
		
		\address{Department of Mathematics,\\ 
			Faculty of Sciences Semlalia,\\
			 Cadi Ayyad University,\\
			  Marrakech, Morocco.\\
			\printead*{e1}}
		
		\address{Department of Mathematics,\\ 
			Faculty of Sciences Semlalia,\\
			Cadi Ayyad University, Marrakech.\\
			Mohammed VI Polytechnic University\\
			Benguerir, Morocco. \\		
			\printead*{e2}\\
			\printead*{e3}
		}
	\end{aug}
\begin{abstract}: In this paper, we study the optimal stopping problem in the case where the reward is given by a family $(\phi(\tau ),\;\;\tau \in \stopo)$ of non negative random variables indexed by predictable stopping times.  We treat the problem by means of Snell's envelope techniques. We prove some properties of the value function family associated to this setting.
\end{abstract}


%\end{keyword}

%\begin{keyword}[class=MSC]
%\kwd[Primary ]{}
%\kwd{}
%\kwd[; secondary ]{}
%\end{keyword}

%\begin{keyword}
%\kwd{}
%\kwd{}
%\end{keyword}
\end{frontmatter}

\section{Introduction}\label{sec1}
\begin{spacing}{1.3}
The classical optimal stopping problem has been studied intensively in many papers in the literature. Without quoting all of them, let us mention the  works of Mertens \cite{Mertens} and \cite{Mertens2}, Bismut and Skalli \cite{Bismut} and  Maingueneau \cite{Main}. For a classical exposition of the optimal stopping problem, we refer to Karatzas and Shreve \cite{KS2} and Peskir and Shiryaev \cite{Peskir}. Generally, the reward of the optimal stopping problem is given by a RCLL process of class $\mathbb{D}$. The most general result in this setting with the right upper semicontinuity as assumption on the reward dates back to  El Karoui \cite{EK}, and has been recently explored to study RBSDEs when the obstacle is not right continuous, in the seminal work \cite{Gro-off} by Grigorova, Imkeller, Offen, Ouknine and  Quenez.  We mention also \cite{Gro-non} in which the authors studied the optimal stopping problem with non linear f-expectation without any regularity assumption on the reward. In all these references cited above, the problem is set and solved in the setup of processes. In \cite{Kob}, Kobylanski and Quenez generalized the classical problem to the case  of a reward family of random variables indexed by stopping times, which is more general than the classical setup of  processes.
Let $S$ be a stopping time,  let
\begin{eqnarray} V(S):=\esssup_{
	\tau\in\mathcal{T}_S}{E[\phi(\tau)|\mathcal{F}_{S^-}]},
\end{eqnarray}
be the value function defined at time $S$ with a reward family $\phi$, where the supremum is taken over the class $\mathcal{T}_S$ of all stopping times $\tau$ such that $ \tau \geq S$ a.s. The objective is to find an optimal stopping time at which  the expected gain reaches its maximum value.
\\
One crucial key in the usual approach consists in the use of the aggregation step of the family $(V(S),\;\;S \in \mathcal{T}_0)$ 
by an optional process $(V_t)$  that is, for each stopping time $S$, $V(S)=V_S$ a.s. The second key consists  in showing that this process is a supermartingale of class $\mathbb{D}$, and thus, it admits a Mertens decomposition which is the analogous of the Doob Meyer decomposition in the right continuous case.
In \cite{Kob}, the approach is purely based on  avoiding aggregation step as well as the use of Mertens decomposition, by considering the setup of family reward, which has appeared as relevant and appropriate as it allows, to release some hypotheses made on the reward.
In \cite{Karoui}, El karoui  gave an extension of \cite{Main} to the predictable case, in which she mentioned the complexity to exhibit conditions on the reward process ensuring the existence of the solution in this framework. 
Inspired from the work of Kobylanski and Quenez \cite{Kob}, we revisit  the general  optimal stopping problem but in the setup of what we call \emph{predictable family} ,  that is a family of random random variables which is indexed by predictable stopping times.
We use an approach which combines some aspects of both approaches. The interest of the predictable general framework has been stressed by Dellacherie in \cite{DelLen2}. Furthermore, this setup establishes a general ground for optimal stopping problems beyond all classical models in continuous time as well as predictable processes setting, which are included as special cases.\\For a predictable stopping time S, the \emph{predictable value family function}   is given by 
\begin{eqnarray}
V_p(S):=\esssup_{
	\tau\in\stops}{E[\phi(\tau)|\mathcal{F}_{S^-}]}.
\end{eqnarray}
where $\stops$ is the set of predictable stopping times $\tau$ with $\tau \geq S$ a.s. 
\\ 
Here   $\phi$ is called  an admissible predictable reward family. In other words, $(\phi(\tau ),\;\;\tau \in \mathcal{T})$ is a family of random variables indexed by stopping times, and verifying two conditions. First, for each  $\tau \in \mathcal{T}_0^p$, $\phi(\tau)$ is an $\mathcal{F}_{\tau ^-}$-measurable $ \bar{\mathbb{R}}^+$-valued random variable. Second, the following condition holds: for each $\tau, \tau'\in \mathcal{T}_0^p $, $\phi(\tau)=\phi(\tau')$ a.s. on $\{\tau=\tau'\}$.
\\
In the first part of the present paper,
we exam some properties of  the family $V_p$. We also prove that this  value family  satisfies an equation in the same spirit of the Bellman equation, from which we derive  many properties as the predictable admissibility, some useful local properties and the supermaringale property of $V_p$. In particular, we characterize the predictable value function  as the predictable Snell envelope system of family $\phi$, defined as the smallest predictable supermartingale family greater than $\phi$. These results highlight the interest of our study of predictable value function, which generalize the notion of strong predictable Snell envelope process. \\
The final section,  describes precisely the difficulties by using  a penalization of Maingueneau \cite{Main}, to obtain $\varepsilon$-optimal predictable stopping times.

\textbf{Notation and terminology:}\\
We start with some notations. We fix a stochastic basis with finite horizon $T \in \mathbb{R}_{+}^*$.
$$(\Omega,\mathbb{F}=(\mathcal{F}_{t})_{t \in {[0,T]}},P).$$
We assume that the filtration $\mathbb{F}$ satisfies the usual assumptions of right continuity and completness. Importantly, we  assume that the filtration is not quasi-left continuous.\\
We suppose that $\mathcal{F}_{0}$ contains only  sets of probability $0$ or $1$. \\
We denote by $\mathcal{T}_0$ the collection of all stopping times  $\tau$ with values in $[0,T]$.\\
We denote by $\stopo$ the collection of all \emph{predictable} stopping times  $\tau$ with values in $[0,T]$.
More generally , we denote $\stops$ (resp. $\mathcal{T}_{S^+}^p$) the class of predictable  stopping times $\tau\in \stopo$
with $S\leq\tau$ a.s.  (resp. $\tau>S $ a.s. on $\{S<T\}$ and $\tau=T$ a.s. on $\{S=T\}$).
\section{Formulation}
\begin{definition}
	A family of random variables $\{\phi(\tau),\; \tau \in \mathcal{T}_0\}$ is said to be a predictable admissible family if it satisfies the following conditions:
	\begin{enumerate}
		\item for all $\tau \in \mathcal{T}_0^p$, $\phi(\tau)$ is an $\mathcal{F}_{\tau ^-}$-measurable $ \bar{\mathbb{R}}^+$-valued random variable,
		\item for all $\tau, \tau'\in \mathcal{T}_0^p $, $\phi(\tau)=\phi(\tau')$ a.s. on $\{\tau=\tau'\}$.
	\end{enumerate}
\end{definition}
\vspace{0.3cm}
In order to simplify, in what follows we use admissible family to mean \emph{predictable} admissible family.\\

In \cite{Karoui}, the reward is given by a predictable process $(\phi_t)$. In this case, the family of random variables defined by $\{\phi(\tau)=\phi_{\tau},\;\tau \in \mathcal{T}_0^p\}$ is admissible.\\
Let $\{\phi(\tau),\; \tau \in \mathcal{T}_0^p\}$ be an admissible family, called reward.
For all $S \in \stopo$,  the value function $V_p$ at time $S$ is defined by:
\begin{eqnarray}\label{value}
V_p(S):=\esssup_{
	\tau\in\stops}{E[\phi(\tau)|\mathcal{F}_{S^-}]}, 
\end{eqnarray}
The strict value function at time $S$ is defined by:
\begin{eqnarray}
\label{strict value} V_p^+(S):=\esssup_{
	\tau\in\stopsp}{E[\phi(\tau)|\mathcal{F}_{S^-}]}, 
\end{eqnarray}
In the interest of kepping this paper self-contained, we prove some results regarding the families  $(V_p(S),\;S \in \mathcal{T}_0^p)$ and  $(V_p^+(S) ,\;S \in \mathcal{T}_0^p)$.
We now state the following proposition:
\begin{Proposition}\label{increas}
	Given any two arbitary predictable stopping times  $S$ and $\theta$ such that  $\theta \in  \mathcal{T}_{S}^p$, the family 	$\{E[\phi(\tau)|\mathcal{F}_{S^-}] \, , \, \tau\in\mathcal{T}_{\theta}^p\;\; (\mbox{resp.}\;\mathcal{T}_{\theta^+}^{p})\}$ is closed under pairwise maximization. Furthermore, 
	there exists a sequence of predictable stopping times $(\tau^n)_{n \in \mathbb{N}}$
	with $\tau^n $ in $ \mathcal{T}_{\theta}^p$ (resp.   $\mathcal{T}_{\theta^+}^{p}$)  such that the sequence $(
	E[\phi(\tau^n)|\mathcal{F}_{S^-}])_{n \in \mathbb{N}}$ converges non-decreasingly to $\displaystyle{\rm ess}\sup_{\tau\in \mathcal{T}_{\theta}^p} E[\phi(\tau)|\mathcal{F}_{S^-}]$ (resp. to  $	\displaystyle{\rm ess}\sup_{\tau\in \mathcal{T}_{\theta^+}^{p}} E[\phi(\tau)|\mathcal{F}_{S^-}]$).	
\end{Proposition}
\dproof
	The arguments are the same for  	$\{E[\phi(\tau)|\mathcal{F}_{S^-}],\tau\in\mathcal{T}_{\theta}^{p}\}$ and $\{E[\phi(\tau)|\mathcal{F}_{S^-}], \tau\in\mathcal{T}_{\theta+}^{p}\}$. We prove the statements only for $\{E[\phi(\tau)|\mathcal{F}_{S^-}] \;,\; \tau\in\mathcal{T}_{\theta+}^{p}\}$. 	For any predictable stopping times $\tau^1$ and $\tau^2$  in $\mathcal{T}_{\theta+}^{p}$,
	write\\ $A:=\{\, E[\phi(\tau^2)|\mathcal{F}_{S^-}]\leq
	E[\phi(\tau^1)|\mathcal{F}_{S^-}]\,\}$ and set
	$$\tau^3:=\tau^1 {\bf 1} _A+\tau^2 {\bf 1} _{A^c}.$$ 
	The fact that  $A \in \mathcal{F}_{S^-} \subset \mathcal{F}_{(\tau^1 \wedge \tau^2)^-}=\mathcal{F}_{(\tau^1)^-} \cap \mathcal{F}_{( \tau^2)^-}$, implies that $A \in\mathcal{F}_{(\tau^1)^- }$ and $A \in\mathcal{F}_{(\tau^2)^-}$. Thus, $\tau^3\in\mathcal{T}_{\theta +}^{p}$,  and by the admissibility of $\phi$, it follows that:\\
	$$ {\bf 1} _AE[\phi(\tau^3)|\mathcal{F}_{S^-}]= 
	E[ {\bf 1} _A \phi(\tau^3)|\mathcal{F}_{S^-}]= E[ {\bf 1} _A \phi(\tau^1)|\mathcal{F}_{S^-}]
	= {\bf 1} _AE[\phi(\tau^1)|\mathcal{F}_{S^-}]\;\;\;\;\mbox{ a.s}.$$
	Similarly, we show that $$ {\bf 1} _{A^c}E[\phi(\tau^3)|\mathcal{F}_{S^-}]={\bf 1} _{A^c}E[\phi(\tau^2)|\mathcal{F}_{S^-}]\;\;\;\; \mbox{ a.s}.$$ 
	Consequently,
	$$E[\phi(\tau^3)|\mathcal{F}_{S^-}]=
	E[\phi(\tau^1)|\mathcal{F}_{S^-}]{\bf 1} _A+
	E[  \phi(\tau^2)|\mathcal{F}_{S^-}]{\bf 1} _{A^c}=
	E[\phi(\tau^1)|\mathcal{F}_{S^-}]\vee
	E[\phi(\tau^2)|\mathcal{F}_{S^-}] \mbox{ a.s},$$
	which shows the stability under pairwise maximization. Thus, by a classical result on essential supremum (see e.g. Neveu \cite{Neveu}), there exists a  sequence of predictable stopping times  $(\tau^n)_{n \in \mathbb{N}} \in\mathcal{T}_{\theta+}^{p} $ such that
$$\displaystyle{\rm ess}\sup_{\tau\in \mathcal{T}_{\theta}^{p+}} E[\phi(\tau)|\mathcal{F}_{S^-}]=\sup_n E[\phi(\tau^n)|\mathcal{F}_{S^-}]\quad  \mbox{ a.s}., $$
	by recurrence, we can define a new sequence of stopping times $(\tilde\tau^n)_{n \in \mathbb{N}} \in\mathcal{T}_{\theta+}^{p}$ by 
	$\tilde\tau^1=\tau^1$, and $\tilde\tau^n$ from $(\tilde\tau^{n-1},\tau^n)$ in the same way as in the definition of $\tau^3$ by $(\tau^{1},\tau^2)$. Hence, we can see that $E[\phi(\tilde\tau^n)|\mathcal{F}_{S^-}]$ converges increasingly to $\displaystyle{\rm ess}\sup_{\tau\in \mathcal{T}_{\theta^+}^p} E[\phi(\tau)|\mathcal{F}_{S^-}]$. 
The proof is thus complete. 
\fproof

As an immediate application of the previous result, for any $S \in \mathcal{T}_0^p$ , $V_p(S)$ ${\mbox{\rm{(resp.}  } } V_p^+(S))$ can be approximated by an appropriate increasing sequence. To be precise, we have the following proposition:
% We formulate this result in the following lemma :
\begin{Proposition}\label{P1.2a}{\em(Optimizing sequences for $V_p$ and $V_p^+$)}
	There exists a sequence of predictable stopping times $(\tau^n)_{n \in \mathbb{N}}$
	with $\tau^n $ in $ \mathcal{T}_{S}^p$ (resp.   $\mathcal{T}_{S^+}^p$),  such that the sequence $(E[\phi(\tau^n)|\mathcal{F}_{S^-}])_{n \in \mathbb{N}}$ is increasing and such that
	
	\[V_p(S) \quad {\mbox{\rm{(resp.}  } } V_p^+(S) \mbox{\rm{)}}\quad = \lim_{n \to \infty} \uparrow E[\phi(\tau^n)|\mathcal{F}_{S^-}] \quad
	\mbox{\rm a.s.}\]
\end{Proposition}
\dproof
	The result follows immediately by taking $\theta=S$ in Proposition \ref{increas}.
\fproof
\begin{lemma}\label{lemm3}
	Let $S \in \stopo $ and $ \theta \in \stops$. Let $\alpha $ be a nonnegative bounded $\mathcal{F}_{\theta^-}$-measurable random variable.  We have,
	\begin{eqnarray}\label{eqBellman}
	E[\alpha V_p(\theta)|\mathcal{F}_{S^-}]=\esssup_{  \tau \in \mathcal{T}_{\theta}^p}E[\alpha \phi(\tau)|\mathcal{F}_{S^-}],
	\end{eqnarray}
	
	and the similar result for the strict value function can re-expressed as the following 
	$$E[\alpha V_p^+(\theta)|\mathcal{F}_{S^-}]=\esssup_{\tau \in \mathcal{T}_{\theta^+}^p}E[\alpha \phi(\tau)|\mathcal{F}_{S^-}].$$
\end{lemma}
\dproof
	Let us prove the result for the value family $V_p$. Let $\tau \in \mathcal{T}_{\theta}^p$,  by iterating expectation and using that $\alpha $ is a nonnegative bounded $\mathcal{F}_{\theta^-}$-measurable random variable,  combined  with  $E[ \phi(\tau)|\mathcal{F}_{\theta^-}] \leq  V_p(\theta)$, we obtain
	$$E[\alpha \phi(\tau)|\mathcal{F}_{S^-}]=E[E[\alpha \phi(\tau)
	|\mathcal{F}_{\theta^-}]|\mathcal{F}_{S^-}]=E[\alpha E[ \phi(\tau)|\mathcal{F}_{\theta^-}]|\mathcal{F}_{S^-}]\leq E[\alpha V_p(\theta)|\mathcal{F}_{S^-}].$$
	By taking the essential supremum over
	$\tau \in \mathcal{T}_{\theta}^p$ in the  inequality, we get 
	$$\esssup_{ \tau \in \mathcal{T}_{\theta}^p}E[\alpha \phi(\tau)|\mathcal{F}_{S^-}]\leq E[\alpha V_p(\theta)|\mathcal{F}_{S^-}].$$
	It remains to prove the reverse inequality $" \leq"$. By Proposition \ref{P1.2a},
	there exists a sequence of predictable stopping times $(\tau^n)_{n\in \mathbb{N}}$
	with 	 $\tau^n $ in $ \mathcal{T}_{\theta}^p$  and such that 
	$$V_p(\theta)=\limn \uparrow E[\phi(\tau^n)|\mathcal{F}_{\theta^-}].$$
	
	Since $\alpha$ is $\mathcal{F}_{\theta^-}-$measurable, we obtain that  $\displaystyle \alpha V_p(\theta)= \limn\uparrow E[ \alpha \phi(\tau^n)|\mathcal{F}_{\theta^-}]$ a.s. Therefore, applying the monotone convergence theorem and the fact that $S \leq \theta$ a.s. we derive that: \\
	$$E[\alpha V_p(\theta)|\mathcal{F}_{S^-}]=\limn \uparrow E[\alpha \phi(\tau^n)|\mathcal{F}_{S^-}].$$
	Hence,
	$$E[\alpha V_p(\theta)|\mathcal{F}_{S^-}] \leq \esssup_{  \tau \in \mathcal{T}_{\theta}^p}E[\alpha \phi(\tau)|\mathcal{F}_{S^-}].$$
	This with the previous inequality leads to the desired result.
\fproof
\begin{remark}
	Note that if $(\phi(\tau),\; \tau \in \stopo)$ is an admissible family,	then for each $S \in \stopo$, and for each $\alpha$ nonnegative bounded $\mathcal{F}_{S^-}$-measurable random variable, the family $(\alpha\phi(\tau),\; \tau \in \stops)$  can be shown $S$-admissible.
\end{remark}
Let $(\phi(\theta),\; \theta \in \stopo)$ be an admissible family. Let $S \in \stopo$, let $\alpha$ be a nonnegative bounded $\mathcal{F}_{S^-}$-measurable random variable. Let $(V^\alpha(\tau),\; \tau \in \stops)$ be the  value function associated with the reward $(\alpha\phi(\theta),\;\theta \in \mathcal{T}_{S}^p)$, defined for each $\tau \in \mathcal{T}_{S}^p $ by 
$$V^\alpha(\tau):=\esssup_{\theta \in \mathcal{T}_{\tau}^p}E[\alpha \phi(\theta)|\mathcal{F}_{\tau^-}].$$
Let $(V^{\alpha+}(\tau),\; \tau \in \stops)$ be the strict  value function associated with the same reward, defined for each $\tau \in \mathcal{T}_{S^+}^p $ by 
$$V^{\alpha+}(\tau):=\esssup_{  \theta \in \mathcal{T}_{\tau+}^p}E[\alpha \phi(\theta)|\mathcal{F}_{\tau^-}].$$
Now, we will state some interesting properties:
\begin{Proposition}\label{Propalpha}
	Let $(\phi(\tau),\; \tau \in \stopo)$ be an admissible family, $S \in \stopo$ and let $\alpha$ be a nonnegative bounded $\mathcal{F}_{S^-}$-measurable random variable. The  value function $(V_p(\tau),\; \tau \in \stops)$  and the strict value function $(V_p^+(\tau),\; \tau \in \stopsp)$ satisfy the following equalities:
	\begin{itemize}
		\item $V^{\alpha}(\tau)=\alpha V_p(\tau)\;\;$ a.s. for all  $\tau \in \stops$.
		\item $V^{\alpha+}(\tau)=\alpha V_p^+(\tau)$  a.s. for all  $\tau \in \stopsp$.
	\end{itemize}
\end{Proposition}
%%%%%%%%%%%%%%%%%%%%%%%%%%%
\dproof
	Let	$\tau \in \stops$ and $\theta \in \mathcal{T}_{\tau}^p$. By the definition of the essential supremum (see Neveu \cite{Neveu} ),	$\alpha E[\phi(\theta)|\mathcal{F}_{\tau^-}]=E[\alpha \phi(\theta)|\mathcal{F}_{\tau^-}] \leq V^{\alpha}(\tau)$.
	Thus, by the characterization of the essential suprmem, we have 
	$\alpha V_p(\tau) \leq  V^{\alpha}(\tau) $. By the same arguments we can show that $V^{\alpha}(\tau) \leq  \alpha V_p(\tau) $. This concludes the proof for the value function $V_p$. The proof is the same for strict value function $V_p^+$.
\fproof
%%%%%%%%%%%%%%%%%%%%%%%%%%%%

Let $S \in \stopo$ and $A \in\mathcal{F}_{S^-}$.  If we take $\alpha={\bf 1}_{A}$, we denote $V^\alpha$ by  $V^A$. Thus, $V^A$  is the value function associated with the  reward $(\phi(\tau)1_A,\; \tau \in \mathcal{T}_{S}^p)$, defined for each $\tau \in \mathcal{T}_{S}^p $ by 
$$V^A(\tau):=\esssup_{  \theta \in \mathcal{T}_{\tau}^p}E[ \phi(\theta)1_A|\mathcal{F}_{\tau^-}].$$
Let $V^{A+}$ be the strict predictable value function associated with the same reward, defined for each $\tau \in \mathcal{T}_{S^+}^p $ by 
$$V^{A+}(\tau):=\esssup_{  \theta \in \mathcal{T}_{\tau^+}^p}E[ \phi(\theta)1_A|\mathcal{F}_{\tau^-}].$$
\begin{lemma}\label{local}
	Let $(\phi(\tau),\; \tau \in \stopo)$ be an admissible family. Let $ \tau,\;\tilde \tau\in \stopo$ and denote $A:=\{\tau=\tilde \tau \}$ and $B:=\{\tau>\tilde \tau \}$. Then
	\begin{itemize}
		
		\item $ V^A(\tau)=V^A(\tilde\tau)\;\;\;\mbox{a.s}\;\;\mbox{and}\;\;\;\;\;V^{A+}(\tau)=V^{A+}(\tilde\tau)\;\;\mbox{a.s}$
	\end{itemize}
	\begin{itemize}
		\item  We have also that:
		\begin{equation} \label{eq.premia}
		E[\phi(\tau)|\mathcal{F}_{\tilde \tau ^-}] \,1_{B} \leq V_p^{B+}(\tilde \tau ) \, \,\,\,\mbox{\rm a.s.}
		\end{equation}
		%Moreover, 
		%Ãƒâ€šÃ‚Âµ\begin{equation}
		%V_p^{B}(\tilde \tau ) \leq V_p^{B+}(\tilde \tau ) \, \,\,\,\mbox{\rm a.s.}
		%\end{equation}
		
	\end{itemize}	
\end{lemma}
\dproof Let us show the result for $V^{A+}$.
	For each $\theta \in \mathcal{T}_{\tau^+}^p$, put $\theta_A=\theta {\bf 1}_A+T{\bf 1}_{A^c}$.
	Since $\tau$ and $\tilde \tau$ are predictable stopping times, we have $A \in \mathcal{F}_{\tau^-} \cap \mathcal{F}_{\tilde\tau^-}$. Thus, $\theta_A $ is predictable, by the admissibility of the family $\phi$, we get:
	\begin{equation*}%\label{eq_toto2}
	\begin{aligned}
	E[\phi(\theta_A){\bf 1}_A| \mathcal{ F}_{\tau^-}]&={\bf
		1}_AE[\phi(\theta)| \mathcal{F}_{\tau^-}]={\bf
		1}_AE[\phi(\theta)| \mathcal{F}_{\tilde\tau^-}]=E[\phi(\theta_A){\bf 1}_A| \mathcal{ F}_{{\tilde\tau}^-}],
	\end{aligned}
	\end{equation*}
	Since $\theta_A \in \mathcal{T}_{\tilde\tau^+} $, we obtain :
	$$E[\phi(\theta_A){\bf 1}_A| \mathcal{ F}_{\tau^-}]\leq  V^{A+}(\tilde\tau).$$
	By arbitrariness of  $\theta  \in \mathcal{T}_{\tau^+}$,  this implies that 
	$$  V^{A+}(\tau)\leq  V^{A+}(\tilde\tau).$$
	By interchanging the roles of $\tau$ and $\tilde{\tau}$, we get $ V^{A+}(\tau)= V^{A+}(\tilde\tau)$.\\
	To prove the second assertion,  let us define the random variable $\overline \tau$  by\\ $\overline \tau := \tau\, 1_{\{\tau >\tilde \tau \}} + T\, 1_{\{\tau \leq \tilde \tau \}}$. Note that $\{\tau >\tilde \tau\} \in \mathcal{F}_{\tilde \tau^-}$, thus $\overline \tau$ belongs to $\mathcal{T}_{\tilde \tau^+}^p$. 
	This combined with the fact that $\{\tau >\tilde \tau\} \in \mathcal{F}_{\tilde \tau^-}$ and the admissibility of the family
	$\phi $ lead to:
	\[E[\phi(\tau)|\mathcal{F}_{\tilde \tau^-}] \,1_{\{\tau >\tilde \tau\}} =E[\phi(\overline \tau)1_{\{\tau >\tilde \tau\}}|\mathcal{F}_{\tilde \tau^-}]=E[\phi(\overline \tau)1_{B}|\mathcal{F}_{\tilde \tau^-}]\leq V_p^{B+}(\tilde \tau) \,\,\,\mbox{\rm a.s.}\]
	Consequently, we get the desired result.
\fproof

Now, we will state te following localization property:
\begin{corollary}\label{corlocal}
	Let $(\phi(\tau),\; \tau \in \stopo)$ be an admissible family, $S \in \stopo$ and let $A \in \mathcal{F}_{S^-}$-measurable random variable. The  value function $(V_p(\tau),\; \tau \in \stops)$  and the strict value function $(V_p^+(\tau),\; \tau \in \stopsp)$ satisfy the following equalities:
	\begin{itemize}
		\item $V^{A}(\tau)={\bf 1}_{A} V_p(\tau)\;\;$ a.s. for all  $\tau \in \stops$.
		\item $V^{A+}(\tau)={\bf 1}_{A} V_p^+(\tau)$  a.s. for all  $\tau \in \stopsp$.
	\end{itemize}
\end{corollary}
\dproof
	The result is a direct application of the Proposition \ref{Propalpha}.
\fproof
\begin{remark}\label{Remark2}
	Let $ \tau,\;\tilde \tau\in \stopo$. Then,
	\begin{equation} \label{eq.p}
	E[\phi(\tau)|\mathcal{F}_{\tilde \tau ^-}] \,1_{\{\tau>\tilde \tau \}} \leq V_p^+(\tilde \tau )\,1_{\{\tau>\tilde \tau \}}\, \,\,\,\mbox{\rm a.s.}
	\end{equation}
\end{remark}
\begin{remark}\label{rem.p}
	Let $S \in \stopo$. Note that if $A \in\mathcal{F}_{S^-}$, we can always decompose the family  $(V_p(\tau),\; \tau \in \stops)$ as the following:
	$$ V_p(\tau)=V^A(\tau)+V^{A^c}(\tau) \;\;\;\;\;\;\; \mbox{for all}\;\;\;\;\;\tau \in \stops.$$ 
\end{remark}
The equalities above are useful, it allows us to prove the admissibility of the value functions $V_p$ and $V_p^+$.
\begin{Proposition}\label{P1.Adm}\emph{(Admissibility of $V_p$ and $V_p^+$)}\\
	The families $V_p=(V_p(S), S\in \stopo)$ and $V_p^+=(V_p^+(S), S\in \stopo)$ defined by (\ref{value}) and (\ref{strict value}) are  admissible.
\end{Proposition}
\dproof Let us show the result for $V_p^+$.
	For each $S\in\stopo, \; V_p^+(S)$ is an $\cf_{S^-}$-measurable random variable, due to the definition of the essential supremum (cf. e.g.  \cite{Neveu}).\\
	Let us prove Property $2$ of the  definition of admissibility. Take $\tau$ and $\tilde \tau$ in $\stopo$. We set $A:=\{\tau=\tilde \tau \}$ and
	we show that $V_p(\tau)=V_p(\tilde \tau)$, $P$-a.s. on $A$. \\
	Thanks to Lemma \ref{local}, $V^{A+}(\tau)=V^{A+}(\tilde\tau)$ a.s.
	Let us remark that $A \in \mathcal{F}_{\tau^- \wedge \tilde\tau^- } $.
	By the second statement of Corollary \ref{corlocal}, we have $$V_p^+(\tau){\bf 1}_{A}=V_p^+(\tilde\tau){\bf 1}_{A} \;\;\;\mbox{a.s.}$$	
	Thus the desired result.
\fproof
\begin{definition}[\emph{Predictable supermartingale system}]
	An admissible family $U:=(U(\tau), \; \tau \in\stopo)$ is
	said to be a  \emph{predictable supermartingale system} (resp. a  \emph{ predictable martingale system}) if, for any 
	$\tau, \tau^{'}$ $ \in$ $\mathcal{T}_0^p$ such that $\tau^{'} \geq \tau$ a.s.,
	\begin{eqnarray*}
		E[U(\tau')|{\cal{F}}_{\tau^-}] \leq  U(\tau) \quad \,\mbox{a.s.}&& {\rm (resp.}, \quad 	E[U(\tau')|{\cal{F}}_{\tau^-}]  =  U(\tau) \quad \,\mbox{a.s.}).
	\end{eqnarray*}
\end{definition} 

\vspace{0.3 cm}
A progressive process $X=(X_t)_{t \in [0,T]}$  is called a predictable strong supermartingale if it is a supermartingale, such that the family $(X_\tau,\; \tau \in \stopo)$ is a predictable supermartingale system.

\begin{corollary}\label{mar.loc}
	Let  $ S \in \stopo$ and $A \in \mathcal{F}_{S^-}$. If the family $(V_p(\tau),\; \tau \in \stops)$ is a predictable martingale system, then the family $(V^A(\tau),\; \tau \in \stops)$ is also  a predictable martingale system.  
\end{corollary}
\dproof
	Let $\tau_1< \tau_2 \in \stops$. Since  $S \leq \tau_1$, we have $A \in \mathcal{F}_{{\tau_1}^-} $. By applying Corollary \ref{corlocal}, and by using the martingale property of the system $(V_p(\tau),\; \tau \in \stops)$, we get  
	$$E[V^A(\tau_2)|\mathcal{F}_{\tau_{1}-}]=E[V_p(\tau_2)1_A|\mathcal{F}_{\tau_1-}]=E[V_p(\tau_2)|\mathcal{F}_{{
			\tau_1-}}]1_A=V_p(\tau_1)1_A=V^A(\tau_1).$$
	This concludes the proof.
\fproof

\begin{lemma}\label{lemm4}
	\begin{itemize}
		\item The admissible families $\{V_p(\tau), \tau \in \stopo\}$ and $\{V_p^+(\tau), \tau \in \mathcal{T}_{0^+}^p\}$ are predictable supermartingale systems. 
		\item The value family $V_p$ is characterized as the  predictable Snell envelope system associated with $\{\phi(S), S\in \stopo\}$, that is, the smallest supermartingale system which is greater (a.s.) than $\{\phi(S), S\in \stopo\}$.
	\end{itemize}
\end{lemma}
\dproof The arguments are the same for $V_p$ and $V_p^+$. Let us prove the first point for $V_p$. Let $S \leq \tau \in \stopo$. 
	Applying Lemma \ref{lemm3}, equation  \ref{eqBellman} holds when $\alpha =1$. Since  $S \leq  \tau$, we get
	$$E[ V_p(\tau)|\mathcal{F}_{S^-}]=\esssup_{ \theta \in \mathcal{T}_{\tau}}E[ \phi(\theta)|\mathcal{F}_{S^-}] \leq \esssup_{ \theta \in \mathcal{T}_S}E[ \phi(\theta)|\mathcal{F}_{S^-}]=V_p(S),$$
	which gives the supermartingale property of $V_p$. \\
	Let us prove the second assertion. Let   $\{V_p'(\tau),\tau\in \stopo\}$  be another  supermartingale system
	such that $V_p'(\tau)\geq \phi(\tau)$ a.s. for all $\tau
	\in\mathcal{T}_S^p$. Thus we have
	$$E[\phi(\tau)|\mathcal{F}_{S^-}]\leq E[V_p'(\tau)|\mathcal{F}_{S^-}] \leq V_p'(S) \quad a.s.$$
	for all $\tau
	\in\mathcal{T}_S^p$. Hence
	by taking the essential supremum over
	$\tau \in\stops$, and by using the definition of $V_p$ we find that
	$$V_p(S)=\esssup_{ \tau \in \stops}E[\phi(\tau)|\mathcal{F}_{S^-}]\leq V_p'(S) \quad a.s. $$ 
	for all $S \in\mathcal{T}_0^p$. This gives the desired result. 
\fproof

%%%%%%%%%%%%%%%%%%%%%%%%%
We state the following property which gives the link between $V_p$, $V_p^+$ and $\phi$.  This corresponds to Proposition D.3 in Karatzas and Shreve \cite{KS2} for right continuous procesesses.
\begin{Proposition}\label{prop.vv+} 
	For all $S\in \mathcal{T}_0^p$, $V_p(S)= \phi(S)\vee V_p^+(S)$ a.s.
\end{Proposition}
\dproof  Note first that $V_p(S)\geq V_p^+(S)$ a.s. and that $V_p(S)\geq \phi(S)$ a.s., which yields the inequality $V_p(S)\geq \phi(S)\vee V_p^+(S)$ a.s. It remains to show the other inequality. 
	Select $\tau$ in $\mathcal{T}_S$. We can rewrite $E[\phi(\tau)|\mathcal{F}_{S^-}]$ as:
	\begin{eqnarray}\label{eq.1}
	E[\phi(\tau)|\mathcal{F}_{S^-}]= E[\phi(\tau)|\mathcal{F}_{S^-}]\, 1 _{\{\tau =S \}}
	+E[\phi(\tau)|\mathcal{F}_{S^-}]\, 1 _{\{\tau >S \}}\,\,\,{\rm a.s.}.
	\end{eqnarray}
	Note that $ 1 _{\{\tau =S \}} \in \mathcal{F}_{S^-}$, thus, by the  admissibility of $\phi$, we get  
	\begin{eqnarray}\label{eq.2}
	E[\phi(\tau)|\mathcal{F}_{S^-}]= \phi(S)\, 1 _{\{\tau =S \}}
	+E[\phi(\tau)|\mathcal{F}_{S^-}]\, 1 _{\{\tau >S \}}\,\,\,{\rm a.s.}.
	\end{eqnarray}
	We have by remark \ref{Remark2}
	\begin{equation} \label{eq.premia}
	E[\phi(\tau)|\mathcal{F}_{S^-}] \,1_{\{\tau >S\}} \leq V_p^+(S) \,1_{\{\tau >S\}} \,\,\,\mbox{\rm a.s.}
	\end{equation}

	Consequently,  by (\ref{eq.2}) and (\ref{eq.premia}) we get:
	$$E[\phi(\tau)|\mathcal{F}_{S^-}]
	\leq \phi(S)\, 1 _{\{\tau =S \}}
	+ V_p^+(S)\, 1 _{\{\tau >S \}}\,\,\,{\rm a.s.}.$$
	\,Therefore,
	\[E[\phi(\tau)|\mathcal{F}_{S^-}]\leq \phi(S)\vee V_p^+(S)\,\,\,{\rm a.s.}\]
	By taking the essential supremum over $\tau \in \mathcal{T}_S^p$, we derive that $V_p(S)\leq \phi(S)\vee V_p^+(S)$ a.s.. The proof is thus complete.
\fproof 

We now state the following  lemma.
\begin{lemma}\label{lemme.un}
	Let $(\phi(\tau),\tau\in \mathcal{T}_0^p)$ be an admissible family.
	For each $\tau, \, S$ $\in$ $\mathcal{T}_0^p$, we have
	\[ E[V_p(\tau)|\mathcal{F}_{S^-}] \leq V_p^+(S)\,\,\,{\rm a.s.} \,\,\, {\rm on} \,\,\,\{\tau > S \}. \]
\end{lemma}
\dproof 
	Select	$\tau \in  \mathcal{T}_0^p$. 
	On account of Proposition \ref{P1.2a}, there exists an optimizing sequence of stopping times $(\tau^n)$ with $\tau^n $ in $ \mathcal{T}_{\tau}^p$  such that 
	$\displaystyle V_p(\tau)  =\lim_{n \to \infty} \uparrow E[\phi( \tau^{n})| \mathcal{F}_{\tau^-}]\quad
	\,\mbox{a.s..}$\\
	Thus, we derive that a.s. on $\{\tau > S \}$, the following equalities hold
	\[ E[V_p(\tau)|\mathcal{F}_{S^-}] = E[\lim_{n \to \infty} \uparrow E[\phi( \tau^{n})| \mathcal{F}_{\tau^-}]|\mathcal{F}_{S^-}]
	=   \lim_{n \to \infty} \uparrow E[\phi( \tau^{n})| \mathcal{F}_{S^-}]\;{\rm a.s.},\]
	we have used here the monotone convergence theorem of conditional expectation.\\  
	Now, on $\{\tau > S \}$,  since
	$\tau^{n} \geq \tau > S$ a.s., in view of remark \ref{Remark2}, we have \\
	$ E[\phi( \tau^{n})| \mathcal{F}_{S^-}] \leq V_p^+(S) $ a.s.
	Passing to the limit in $n$ and using the previous equality, we obtain that $E[V_p(\tau)|\mathcal{F}_{S^-}]  \leq V_p^+(S)$ a.s. on $\{\tau > S \}$.
\fproof  

\section{Regularity and predictable value function }
%\begin{definition}
%	An admissible family $(\phi(\theta),\theta \in \stopo)$ is said to be right continuous along stopping times %$(\theta_n)_{n \in \mathbb{N}}$ such that $\theta_n \downarrow \theta$ one as $\phi(\theta)= \lim_{n \rightarrow %} \phi(\theta_n)$ a.s.	
%\end{definition}
In this paragraph,  we focus our attention to provide many useful properties for the  value and strict value functions in the predictable setting. First, let us introduce these new definitions:
\begin{definition}
	An admissible family $(\phi(\theta),\;\theta \in \mathcal{T}_0^p)$ is said to be right continuous along  predictable stopping times (RCP) if for any $\theta \in \stopo$ and for any sequence of predictable stopping times $(\theta_n)_{n \in \mathbb{N}}$ such that $\theta_n \downarrow \theta$ one has $\phi(\theta)= \limn \phi(\theta_n)$.
\end{definition}
\begin{definition}
	An admissible family $(\phi(\theta),\;\theta \in \mathcal{T}_0^p)$ is said to be right continuous along  predictable stopping times  in expectation (RCPE) if for any $\theta \in \stopo$ and for any sequence of predictable stopping times $(\theta_n)_{n \in \mathbb{N}}$ such that $\theta_n \downarrow \theta$ one has $E[\phi(\theta)]= \limn E[\phi(\theta_n)]$.
\end{definition}
\begin{definition} \label{defr} An admissible family $(\phi(\theta),\;\theta \in \mathcal{T}_0^p)$ is said to be {\em left-upper semicontinuous (l.u.s.c.) along stopping times}  if for all $\theta \in {\cal T}_0$ and for each non decreasing sequence of stopping times $ (\theta_n)$ such that $\theta^n \uparrow \theta$ a.s.\,,
	\begin{equation}\label{usc}
	\phi(\theta) \geq \limsup_{n\to \infty} \phi(\theta_n).
	\end{equation}
\end{definition}
\begin{definition}
	Let $S$ $\in$ $\mathcal{T}_0$. An admissible family $(\phi(\tau),\tau\in \mathcal{T}_0)$ is said to be  \emph{right limited along stopping times} (RL) at $S$ 
	if 
	there exists an $\mathcal{F}_{S}$-measurable random variable $\phi(S^+)$ such that,
	for any non increasing sequence  of stopping times  $(S_n)_{n\in\mathbb{N}}$,
	such that  $S_n\downarrow S$ and $S_n >S$ for each $n$, one has
	$\displaystyle{\phi(S^+)=\lim_{n \to\infty} \phi(S_n)} $.
\end{definition}
The following result of Kobylanski et al. see Theorem $4.6$ \cite{Kob}, can be expressed in our setting by the following:
\begin{theorem} \label{surl}
	A   predictable supermartingale family $(U(\tau), \tau \in \mathcal{T}_0)$, with $U(0) < + \infty$, is right limited along stopping times (RL) and left limited  along stopping times (LL)
	at any   stopping time $S$. 
\end{theorem} 
\begin{Proposition}\label{prop3}
	Let $(\phi(\tau),\tau \in \mathcal{T}_0^p)$ be an uniformly integrable admissible family.
	For each $\tau, \, S$ $\in$ $\mathcal{T}_0^p$ such that $\tau >S$, one has 
	\begin{itemize}
		\item  $E[V_p(\tau)|\mathcal{F_{S^-}}] \leq E[V_p(S^+)|\mathcal{F_{S^-}}] $ \;{\rm a.s.}
		\end{itemize}
	Moreover, we have
		\begin{itemize}
		\item $V_p^+(S) \leq E[V_p(S^+)|\mathcal{F_{S^-}}]$\;{\rm a.s.}
		\item $ V_p(S) \leq V_p(S^-) $ \;{\rm a.s.}
	\end{itemize}
\end{Proposition}
	\dproof
		Let $ \tau >S \in \stopo$ and  $(S_n:=S+\frac{1}{n} \wedge \tau)$, $S_n$ is a predictable stopping time for all $n \in  \mathbb{N}$. By the supermartingale property, the RL property and the uniform integrability of the family $(V_p(S_n))$, we obtain
		\begin{eqnarray*}
			E[V_p(\tau)|\mathcal{F_{S^-}}] 
			&\leq& \limn E\left[E\left[V_p(\tau)| \mathcal{F}_{S_n ^-}\right]|\mathcal{F_{S^-}}\right]\\
			&\leq& \limn E\left[ V_p\left(S_n \right)|\mathcal{F_{S^-}}\right] \\
			&\leq & E[V_p(S^+)|\mathcal{F_{S^-}}] \;{\rm a.s.}.
		\end{eqnarray*}
		For the second assertion, note that $\phi(\tau)\leq V_p(\tau)$ a.s. Thus,
		$E[\phi(\tau)|\mathcal{F_{S^-}}] \leq E[V_p(\tau)|\mathcal{F_{S^-}}]\;$ a.s.
		This yields by the first assertion, 
		$E[\phi(\tau)|\mathcal{F_{S^-}}] \leq E[V_p(S^+)|\mathcal{F_{S^-}}]$ a.s. By arbitrariness of $ \tau $ and since $\tau>S$ we get,
		$V_p^+(S) \leq E[V_p(S^+)|\mathcal{F_{S^-}}]$ a.s.\\
		Now, let us prove the third assertion. 	Let $(S_n)$ be a sequence of predictable stopping times foretelling $S$. By the predictable supermartingale property of the family $(V_p(\tau),\;\tau \in \stopo)$, $E[V_p(S)|\mathcal{F}_{S_n^-}]\leq V_p(S_n)$ a.s. By letting $n$ tend to $\infty$, we obtain $E[V_p(S)|\vee\mathcal{F}_{S_n^-}]\leq\lim_n V_p(S_n)$ a.s. On the other hand, $\mathcal{F}_{S^-}=\vee \mathcal{F}_{S_n^-}$. Since the family $V_p$ is admissible and uniformly integrable, we get from																																																																																										 the last inequality that $V_p(S)\leq V_p(S^-)$ a.s., which proves the third statement.
	\fproof

\begin{Proposition}\label{propV_p}
	Let $(\phi(\tau),\tau \in \mathcal{T}_0^p)$ be an uniformly integrable admissible family. For each $S \in \stopo$
	\begin{itemize}
		\item $E[V_p(S^+)|\mathcal{F}_{S^-}]\leq V_p(S)$ a.s.
		\item If $(V_p(\tau),\;\tau \in \stopo)$ is RCPE, then, 
		$$\forall S \in \stopo,\;E[V_p(S^+)|\mathcal{F}_{S^-}]=V_p(S)\;\;\;\;\;\mbox{a.s.}$$
	\end{itemize}
\end{Proposition}
\dproof 
	Let $S \in \mathcal{T}^p_0$, let $(S_{n})_{n} \in \mathcal{T}^p_{S^+}$ such that $S_n \downarrow S$. By the predictable supermartingale property of the family $(V_p(\tau),\;\tau \in \stopo)$, $E[V_p(S^n)|\mathcal{F}_{S^-}]\leq V_p(S)$ a.s. Since the family $V_p$ is uniformly integrable, we get by passing to the limit in the last inequality that $E[V_p(S^+)|\mathcal{F}_{S^-}]\leq V_p(S)$ a.s., which proves the first statement.\\
	Now, let us prove the second statement. \\
	Since by assumption  $V_p$ is RCPE, $E[V_p(S)]=\limn E[V_p(S^n)]$.  Furthermore, $V_p$ is right limited at any predictable stopping time, thus the uniform integrability property yields,   $E[V_p(S^+)]=\limn E[V_p(S^n)]$. Therefore, $E[V_p(S^+)]=E[E[V_p(S^+)|\mathcal{F}_{S^-}]]=E[V_p(S)]$, this combined with this first inequality $E[V_p(S^+)|\mathcal{F}_{S^-}]\leq V_p(S)$ a.s., give the second point.
\fproof
\begin{corollary}\label{corol3}
	Let $(\phi(\tau),\tau \in \mathcal{T}_0^p)$ be an uniformly integrable admissible family.
	For each $S \in \stopo$
	We have 
	\begin{itemize}
		\item $V_p(S)=\phi(S)\vee E[V_p(S^+)|\mathcal{F_{S^-}}]$ a.s.
		\item $V_p(S)-E[V_p(S^+)|\mathcal{F_{S^-}}]=(
		V_p(S)-E[V_p(S^+)|\mathcal{F_{S^-}}]){\bf 1}_{\{V_p(S)=\phi(S)\}}$ a.s. 
	\end{itemize}		
\end{corollary}
\dproof
	For all $S \in \stopo$, $\phi(S) \leq V_p(S)$ a.s and by Proposition \ref{propV_p},  $E[V_p(S^+)|\mathcal{F}_{S^-}] \leq V_p(S)$ a.s. 
	Applying the Proposition \ref{prop.vv+} and the second assertion of Proposition  \ref{prop3},  we get 
	$$ V_p(S)=\phi(S)\vee V_p^+(S) \leq \phi(S)\vee E[V_p(S^+)|\mathcal{F_{S^-}}] \leq V_p(S)\;\mbox{a.s}.$$
	
	Thus, the inequalities become equalities:
	$$V_p(S)=\phi(S)\vee E[V_p(S^+)|\mathcal{F_{S^-}}] \;\mbox{a.s}  .$$
	The second statement is a direct consequence of the first assertion.
\fproof
\begin{Proposition}\label{prop.vv-8} 
	Let $(\phi(\tau),\tau \in \mathcal{T}_0^p)$ be a left limited admissible family such that $\sup_{\tau \in \stopo}E[\phi(\tau)]< \infty$. We have	for all $\tau\in \mathcal{T}_
	{0^+}^p$, $V_p(\tau^-)= \phi(\tau^-)\vee V_p(\tau)$ a.s.
\end{Proposition}
\dproof Let $\tau$ be a predictable stopping time.  Note first that by the third assertion of Proposition \ref{prop3},  $V_p(\tau^-)\geq V_p(\tau)$ a.s. and that $V_p(\tau^-)\geq \phi(\tau^-)$ a.s., we obtain the inequality $V_p(\tau^-)\geq \phi(\tau^-)\vee V_p(\tau)$ a.s. It remains to show the other inequality. Let $(\tau^n)$ be a sequence of predictable stopping times foretelling $\tau$. Then
\begin{eqnarray*}
& E[V_p(\tau^-)]=\lim_nE[V_p(\tau^n)]\leq  \lim_{S \geq\tau_n }E[\phi(S)\, 1 _{\{\tau >S \geq\tau^n \}}
+\phi(S)\, 1 _{\{\tau \leq S \}}]\\
&\leq  \lim_{S \geq\tau_n }E[\phi(S)\, 1 _{\{\tau >S \geq\tau^n \}}
+V(S)\, 1 _{\{\tau \leq S \}}]\\
& \leq \lim_{\tau >S \geq\tau_n }E[\phi(S) \vee V(\tau)]=E[\phi(\tau^-) \vee V(\tau)].
\end{eqnarray*}
Thus the desired result.
\fproof 
\begin{remark}
	If $\phi$ is an admissible family, then $V_p(S)=E[\phi(S)| \mathcal{F}_{S^-}]\vee E[V_p(S^+)|\mathcal{F_{S^-}}]$ a.s.
\end{remark}
\begin{Proposition}\label{prop-process}
If the family $(V_p(\tau),\;\tau \in \mathcal{T}^p_0)$ is uniformly integrable.																																									Then the family $(E[V_p(\tau^+)|\mathcal{F}_{\tau^-}],\;\tau \in \mathcal{T}^p_0)$ is a strong predictable supermartingale system.
\end{Proposition}
\dproof Let  $\psi(\tau):=E[V_p(\tau^+)|\mathcal{F_{\tau^-}}]$  for all $\tau \in \mathcal{T}_0^p$. Let $\tau,S \in \mathcal{T}^p_0 $ such that $\tau > S$ a.s. By using the definition of $\psi$ and by iterating expectation,  and that $E[V_p(\tau^+)|\mathcal{F_{\tau^-}}] \leq V_p(\tau)$ a.s, we get
	\begin{eqnarray}
	E[\psi(\tau)|\mathcal{F}_{S^-} ]=E[E[V_p(\tau^+)|\mathcal{F_{\tau^-}}]|\mathcal{F}_{S^-}]\leq E[  V_{p}(\tau)|\mathcal{F}_{S^-}]\;\mbox{a.s}.
	\end{eqnarray}
	By the first assertion in Proposition \ref{prop3}, $E[  V_{p}(\tau)|\mathcal{F}_{S^-}]\leq E[V_p(S^+)|\mathcal{F_{S^-}}]$ a.s., thus, $E[\psi(\tau)|\mathcal{F}_{S^-} ] \leq \psi(S) .$
\fproof

In the follwing, we use some of the results provided above to derive some fine results in the setting of processes.
\begin{Proposition}\label{prop-process}
	Let $(U_{\theta},\;\theta \in \mathcal{T}^p_0)$ be an uniformly integrable predictable value  function. Then 
	\begin{enumerate}
		\item The family $({}^p U_{\theta}^+,\;\theta \in \mathcal{T}^p_0)$ is a strong predictable supermartingale system.
		\item   For each $S \in \mathcal{T}^p_0\;$, ${}^pU_{S}^+\leq U_{S}$ a.s.
		\item If $(U_{\theta},\;\theta \in \mathcal{T}_0^p)$ is RCPE, then, $ {}^p U^+\equiv U.$
	\end{enumerate}
\end{Proposition}
\dproof
First, Let us denote $Z$ the process defined by $Z_\theta:={}^pU_{\theta}^+$  for all $\theta \in \mathcal{T}_0^p$. By using the definition of $Z$ and by iterating expectation and applying the definition of the predictable projection combined with the first assertion of the Corollary \ref{corol3}, we get for $\theta \in \mathcal{T}_{S^+}^p$
	\begin{eqnarray}
	E[Z_{\theta}|\mathcal{F}_{S^-} ]=E[E[U_{\theta^+}|\mathcal{F}_{\theta^-}]|\mathcal{F}_{S^-}]\leq E[U_{\theta}|\mathcal{F}_{S^-}].
	\end{eqnarray}
	By the first assertion in Proposition \ref{prop3}, we have $E[  U_{\theta}|\mathcal{F}_{S^-}]\leq E[U_{S^+}|\mathcal{F_{S^-}}]$. Thus, 
	$$E[Z_{\theta}|\mathcal{F}_{S^-}] \leq E[U_{S^+}|\mathcal{F_{S^-}}]=Z_S.$$

	Thanks to Proposition \ref{propV_p} applied to the predictable supermartingale family $(U_{\theta},\;\theta \in \mathcal{T}^p_0)$ combined with the section theorem , we get the second and the third statements.
\fproof
\section{Optimality criterion and Snell envelope system}\label{section4}
\begin{definition}
	A predictable stopping time $\tau^*  \in \stops$ such that $E[\phi(\tau^*)]< \infty$ is said to be optimal for $V_p(S)$ if and only if 
		\[E[V_p(S)]=\sup_{\tau \stops}E[\phi(\tau)]= E[\phi(\tau^*)] .\]
\end{definition}
We now in position to provide necessary and sufficient conditions, for predictable optimal stopping time, in terms of appropriate martingales. This represents the predictable analogous of Bellman optimality criterium  (c.f El Karoui \cite{Karoui} in the setup of processes).
\begin{Proposition}\label{prop.criterion}\emph{(Optimality criterion) }
	Let $S$ $\in$ $\mathcal{T}_0^p$ and let $\tau_{*} \in \mathcal{T}_S^p$. 	$\tau_{*}$ is $S$-optimal for $V_p(S)$ if and only if
	the following assertions hold:
	\begin{enumerate}
		
		\item $V_p ( \tau_{*}) = \phi( \tau_{*}).$ a.s.
		\item  The family $(V_p(\tau),\; \tau \in \mathcal{T}_{[S, \tau^* ]})$ is a predictable martingale system. 
	\end{enumerate}
\end{Proposition}
\dproof 
	By definition, $\tau^*$ is optimal if and only if
		\[E[V_p(S)]=\sup_{\tau \stops}E[\phi(\tau)]= E[\phi(\tau^*)]. \]
	 Since the value function $V_p$ is a strong predictable supermartingale family greater that $\phi$, we have clearly 
	\[E[V_p(S)]= E[\phi(\tau^*)] = E[V_p(\tau^*)]=E[V_p(\tau^*\wedge \theta)] \;\mbox{for all  }\; \theta \in \stopo \]
	These equalities are equivalent to the conditions of the theorem.
\fproof

Let $(\phi(\tau),\tau\in \mathcal{T}_0^p)$ be an uniformly integrable admissible family. For each $S \in \mathcal{T}_0^p$, suppose that $\hat \tau$ is a predictable optimal stopping time for $V_p(S)$, then, as a consequence of the optimality criterion , the family $(V_p(\tau),\; \tau \in \mathcal{T}_{[S,\hat \tau ]})$ is a	 predictable martingale family. Consider the set
$$\mathcal{A}^p_S=\{\tau \in \mathcal{T}_S^p, \;\;\mbox{such that}\;\;(V_p(\tau), \tau \in \mathcal{T}_{[S,\tau]})\;\;\mbox{is a predictable martingale family}\}.$$ 
\begin{lemma}\label{lem-5}
	For each $S \in \stopo$, the set $\mathcal{A}^p_S$ is stable by pairwise maximization.
\end{lemma}
\dproof
	Let $S \in \stopo$ and $\tau_1,\;\tau_2 \in \mathcal{A}_S^p$. First, we have $\tau_1 \vee \tau_2  \in \mathcal{T}_S^p $. Let us show that  $\tau_1 \vee \tau_2  \in \mathcal{A}^p_S $. This is equivalent to show that $(V_p(\tau), \tau \in \mathcal{T}_{[S,\tau_1 \vee \tau_2]})$  is a predictable martingale family. We have a.s.
	\begin{eqnarray}\label{eqmar}
	E[V_p(\tau_1 \vee \tau_2)|\mathcal{F}_{S^-}]=E[V_p(\tau_2)1_{\{\tau_2>\tau_1\}}|\mathcal{F}_{S^-}]+ E[V_p(\tau_1)1_{\{\tau_1\geq \tau_2\}}|\mathcal{F}_{S^-}].	
	\end{eqnarray}	
	Let $A=1_{\{\tau_2>\tau_1\}}$, thus the equality \ref{eqmar} can be rewritten as: 
	\begin{eqnarray}\label{eqgen}
	E[V_p(\tau_1 \vee \tau_2)|\mathcal{F}_{S^-}]=E[V^A(\tau_2)|\mathcal{F}_{S^-}]+ E[V^{A^c}(\tau_1)|\mathcal{F}_{S^-}].	
	\end{eqnarray}
	
	Since $\tau_2 \in \mathcal{A}_S$ and $A \in \mathcal{F}_{(\tau_1 \wedge \tau_2)^- }$, we have by Corollary \ref{mar.loc}, $(V^A(\tau), \tau \in \mathcal{T}_{[\tau_1\wedge \tau_2 ,\tau_2]})$ ia predictable martingale family.
	Therefore, by iterating expectation, and using that $\tau_1 \in \mathcal{A}_S^p$ and Corollary \ref{corlocal}
	$$E[V^A(\tau_2)|\mathcal{F}_{S^-}]=E[E[V^A(\tau_2)|\mathcal{F}_{(\tau_1 \wedge \tau_2)^-}]|\mathcal{F}_{S^-}]=E[V^A((\tau_1 \wedge \tau_2))|\mathcal{F}_{S^-}]=E[V^A(\tau_1)|\mathcal{F}_{S^-}]. $$
	Hence, the equality \ref{eqgen}, can be expressed as 
	\begin{eqnarray}
	E[V_p(\tau_1 \vee \tau_2)|\mathcal{F}_{S^-}]=E[V^A(\tau_1)|\mathcal{F}_{S^-}]+ E[V^{A^c}(\tau_1)|\mathcal{F}_{S^-}].	
	\end{eqnarray}
	By Remark \ref{rem.p} and using that $\tau_1 \in \mathcal{A}_S^p$, we get
	\begin{eqnarray}
	E[V(\tau_1 \vee \tau_2)|\mathcal{F}_{S^-}]=V^A(S)+V^{A^c}(S)=V(S)\;\;\;\;\mbox{a.s.}
	\end{eqnarray}	
	This yields the desired result.
\fproof

Let us consider the random variable $ \tilde \tau(S) $  defined by
$$ \tilde \tau(S):=\esssup \mathcal{A}^p_S.$$
\begin{lemma}
The random variable	$ \tilde \tau(S)$ is a predictable stopping time.  Moreover, 
	assume that  $(\phi(\tau),\tau\in \mathcal{T}_0^p)$ is  l.u.s.c along stopping times. Then, the family $V_p$ is a predictable martingale system on $[S,\tilde\tau(S)]$	.
\end{lemma}
\dproof
By Lemma \ref{lem-5}, there exists a sequence $\tau^n$ such that $\tau^n \uparrow \tilde  \tau(S)$. Since the family $\phi$ is l.u.s.c, we get from Proposition \ref{prop.vv-8}, that the family $V_p$ is left continuous at $\tilde\tau(S)$. Thanks to the uniform integrability of $(V_p(\tau^n))$, we have
$$E[V_p(\tilde\tau(S))]=\limn E[V_p(\tau_n)]=E[V_p(S)].$$ 
Where the last equality follows from the fact that $(\tau^n)_{n \in \mathbb{N}} \in \mathcal{A}^p_S. $ This concludes the proof.
\fproof

Fix $S \in \mathcal{T}_0^p$.
In order to tackle the existence of optimal stopping times for the value function in the classical case, we proceed to construct a family of "approximately optimal" stopping times. Let $S \in \mathcal{T}_0^p$. For $\alpha \in ]0,1]$, let us introduce the following $\mathcal{F}_{S}$-measurable random variable 
%\begin{equation}\label{tlpS} \tau^{\alpha}(S):=
%{\rm ess} \sup \;\mathbb{T}^{\alpha}_S \quad \mbox{where} \quad \mathbb{T}^{\alpha}_S :=\{\,\tau \in \mathcal{T}^p_{[S,\tilde \tau(S)]}\, , \,\phi(\tau^-)<\alpha %V(\tau^-)  \,\,\mbox {a.s.} \,\}.
%\end{equation}
\begin{equation}\label{tlpS} \tau^{\alpha}(S):= \essinf \;\mathbb{T}^{\alpha}_S \quad \mbox{where} \quad \mathbb{T}^{\alpha}_S :=\{\,\tau \in \mathcal{T}_S^p\, , \,\alpha V_p(\tau)\leq \phi(\tau)  \,\,\mbox {a.s.} \,\}.
\end{equation}
Since the  set $\mathbb{T}^{\alpha}_S$ is stable by pairwise minimization, for every  $\alpha \in ]0,1[$, $\tau^{\alpha}(S)$ can be approached by a decreasing sequence in this family. As a result,  $\tau^{\alpha}(S)$ defines a stopping time in $\mathcal{T}_S$.
Now we will recall the following definition of right upper semicontinuity along stopping times in expectation:
\begin{definition}
	An admissible family $(\phi(\theta),\;\theta \in \mathcal{T}_0)$ is said to be right upper semicontinuous along stopping times in expectation (right USCE), if for all sequence of stopping times $(\theta_n)_{n \in \mathbb{N}}$ such that $\theta_n \downarrow \theta$ one has $$\limsup_{n\to \infty} E[\phi(\theta_n)]\leq E[\phi(\theta)].$$
\end{definition}
\begin{lemma} \label{stepun}
	Suppose  the reward family $( \phi(\tau), \tau\in \mathcal{T}_0^p)$ is right USCE and $V_p(0)< \infty$. Then, for each $\alpha \in ]0,1[$:
	\begin{equation*}
	\alpha V_p (\tau^{\alpha}(S) )\leq  \phi(\tau^{\alpha}(S) ) \,\,\,{\rm a.s.}
	\end{equation*}
\end{lemma}

Before giving the proof of this lemma, we will recall the following theorem (c.f. Dellacherie Theorem $15$ \cite{DelLen2} ):
\begin{theorem} \label{thmfso1}
	Let 
	$\{\phi(\tau),\;\tau \in \mathcal{T}_0^p\}$ be an admissible family of random variables. Suppose that $\sup_{\tau \in \stopo}E[\phi(\tau)]< \infty$.
	Then there exists a unique strong predictable  supermartingale
	denoted by $(V_t)$ which aggregates the value family $\{V_p(S), S\in \stopo\}$, that is, for each stopping time $S$, $V_p(S)=V_S$ a.s.  Moreover, $V$ is the strong predictable Snell envelope satisfying $V \geq \phi$, i.e. if $V' \geq\phi$ is another strong predictable supermartingale process, then $V \leq V'$.
\end{theorem}
Proof of Lemma \ref{stepun}:
Fix $S \in \stopo$. To simplify the notation, in the following, the stopping time $\tau^{\alpha}(S)$ will be denoted b $\tau^{\alpha}$. The families $(V_p(\tau), \tau\in \mathcal{T}_0^p)$ and $(E[V_p(\tau^+)|\mathcal{F}_{\tau^-}], \tau\in \mathcal{T}_0^p)$ are predictable supermartingales system. It follows  by  Theorem $15$ in \cite{DelLen2}, that these families can be aggregated by predictable processes. Let us denote $V$  the predictable process which aggregates the family $(V_p(\tau), \tau\in \mathcal{T}_0^p)$. We have by Corollary \ref{corol3} combined with section theorem c.f. \cite{DM1}, that:
$$V_{\tau^{\alpha}}-E[V_{(\tau^{\alpha})^+}|\mathcal{F}_{\tau^{\alpha^-}}]=(
V_{\tau^{\alpha}}-E[V_{(\tau^{\alpha})^+}|\mathcal{F}_{\tau^{\alpha^-}}]){\bf 1}_{\{V_p(\tau^{\alpha})=\phi(\tau^{\alpha})\}} a.s.$$
Using again the aggregation equality, we get $$V_p(\tau^{\alpha})-E[V_p((\tau^{\alpha})^+)|\mathcal{F}_{\tau^{\alpha^-}}]=(
V_p(\tau^{\alpha})-E[V_p((\tau^{\alpha})^+)|\mathcal{F}_{\tau^{\alpha^-}}]){\bf 1}_{\{V_p(\tau^{\alpha})=\phi(\tau^{\alpha})\}}a.s.$$
Let $A \in {\cal{F}}_{\tau^{\alpha-}}$, we obtain $$E(\alpha\;V_p(\tau^\alpha)1_{A} )= E( \alpha\; V_p((\tau^{\alpha})^+)1_{A \cap\{ V_p(\tau^{\alpha})> \phi(\tau^{\alpha})\}})+E(\alpha V_p((\tau^{\alpha}))1_{A \cap\{V_p((\tau^{\alpha})=  \!\phi(\tau^{\alpha})\}}).$$
By definition of $\tau^{\alpha}$, there exists a non-increasing sequence $(\tau^n)$  in $\mathbb{T}^{\alpha}_S$ verifying  $\displaystyle{\tau^\alpha= \lim_{n\to \infty} \downarrow \tau^n}$  such that, we have  for each $n$,
\begin{equation}\label{eq.vn}
\alpha V_p(\tau^n) \leq  \!\phi(\tau^n)\,\,\,{\rm a.s.}
\end{equation}
On the other hand, the family $(V_p(\tau^n))_{n \in \mathbb{N}}$ is uniformly integrable, thus
$$ E[1_{A \cap \{V_p(\tau^\alpha)>\phi(\tau^\alpha)\}}V_p((\tau^{\alpha})^+)]= \lim_{n\to \infty} E[1_{A \cap \{V_p(\tau^\alpha)> \phi(\tau^\alpha)\}}V_p(\tau^n) ].$$ 
Hence,
\begin{eqnarray*}
	E(\alpha\;V_p(\tau^\alpha)1_{A} )
	&\leq& \limsupn E(\phi(\tau^n)1_{A \cap\{V_p((\tau^\alpha)> \!\phi(\tau^\alpha)\}})+E( \phi(\tau^{\alpha})1_{A \cap\{V_p(\tau^\alpha)=\!\phi(\tau^\alpha)\}})\\
	&\leq& \limsupn E( \phi(\bar \tau^n)1_A)
\end{eqnarray*}
where $\bar{\tau}^n:=\tau^n 1_{A \cap\{ V_p(\tau^\alpha)> \!\phi(\tau^\alpha)\}}+\tau^\alpha 1_{A \cap\{V_p(\tau^\alpha)=\!\phi(\tau^\alpha)\}} $.  Note that $\bar\tau^n$ is a non-increasing sequence of stopping times which verifies $\tau^{\alpha}=\displaystyle\lim_{n\to \infty} \downarrow \bar \tau^n$.
It follows by the right USCE assumption on the reward family $\phi$ that:
$$E(\alpha\;V_p(\tau^\alpha)1_{A} ) \leq  limsup_{n\to \infty}E[\phi(\bar\tau^n)1_A]\leq E[\phi(\tau^\alpha)1_A].$$	
This holds for each $A \in {\cal{F}}_{\tau^{\alpha-}}$. Hence the desired result.
\fproof
%%%%%%%%%%%%%%%%%MERTENS%%%%%%

In \cite{Mertens}, Mertens gives the analogous of the Doob-Meyer decomposition theorem in the general case. In the following, we will recall the  so-called  \emph{Mertens decomposition} (see Meyer \cite{Me}).

\begin{theorem}[Mertens decomposition]\label{thm_Mertens_decomposition}
	Let $X$ be a strong predictable supermartingale of class $(\mathcal{D})$. 
	There exists a unique uniformly integrable martingale $(M_t)$, a unique predictable right-continuous nondecreasing process $(A_t)$ with $A_0=0$ and $E[A_T] < \infty$, and a unique right-continuous adapted nondecreasing process $(C_t)$, which is purely discontinuous, with $C_{0-}=0$ and $E[C_T] < \infty$, such that 
	\begin{equation}\label{eq_thm_Mertens_decomposition}
	X_t=M_{t^-}-A_{t}-C_{t-}, \,\,\, 0 \leq t\leq T  \,\,\, {\rm a.s.}\\
	\end{equation}
	and 
	\begin{equation}\label{eq_thm_Mertens_}
	\Delta C_t ={} ^{p}X_{t}^+-X_t \,\,\, {\rm a.s.}  \,\,\, \mbox{and}  \,\,\, \Delta A_t= X_{t^-}-X_t.
	\end{equation}
\end{theorem}
\begin{lemma}\label{lemeq}
	Let $(\phi(\tau), \tau \in \stopo)$ be an admissible family with $V_p(0)< \infty$. Let $V$ be the predictable process which aggregates the family $V_p$, having the  following Mertens decomposition:
	\begin{equation*}
	V_t = M_{t^-} - A_{t} - C_{t-}, \,\,\, 0 \leq t \leq T  \,\,\, {\rm a.s.}\\
	\end{equation*}
	
	Then, for each $\alpha  \in ]0,1[$ and for each $S \in \stopo$, 
	\begin{equation}\label{eq.mar}
	V_p(S)=E[M_{\tau^\alpha(S)}-A_{\tau^\alpha(S)}-C_{\tau^\alpha(S)^{-}}|\mathcal{F}_{S^-}]\;\;\;\; \mbox{a.s.}
	\end{equation}
	Moreover,
	$A_{\hat \tau(S)}=A_S$ and  $C_{\hat \tau(S)^-}=C_{S^-}$.
\end{lemma}

\dproof
	Let $\hat V$ be the  Snell envelope family of the family $V_p {\bf 1}_{\{\ \alpha V_p \leq \phi\}}$.\\
	We will show that $V_p(S)=\hat V(S)$ for all $S \in \stopo$. Let  $S \in \stopo$, we can rewrite $\phi(S)$ as :
	\begin{eqnarray*}
		\phi(S)&=& \phi(S){\bf 1}_{\{\phi(S)< \alpha V_p(S)\}}+{\bf 1}_{\{\alpha V_p(S) \leq \phi(S)\}}
		\alpha \phi(S) +(1-\alpha) \phi(S){\bf 1}_{\{\alpha V_p(S) \leq \phi(S)\}} \;\mbox{a.s.}
\end{eqnarray*}
By using the definition of $V$ and $\hat V$, we get
\begin{eqnarray*}
	\phi(S)	&\leq & \alpha V_p(S)+ (1-\alpha) \hat V(S) {\bf 1}_{\{\alpha V_p(S) \leq \phi(S)\}}+ (1-\alpha) \hat V(S) {\bf 1}_{\{\phi(S) < \alpha V_p(S)
			\}}\mbox{a.s.}.
	\end{eqnarray*}
	This yields , 
	$$\phi(S) \leq\alpha V_p(S)+ (1-\alpha) \hat V(S).$$
	Since $ (\alpha V_p(S)+ (1-\alpha) \hat V(S);\;\;S\in \mathcal{T}_{0}^p )$,  is a supermartingale family which is greater than $\phi$, we get
	$$V_p(S) \leq \alpha V_p(S)+ (1-\alpha) \hat V(S).$$
	Hence, we obtain $V_p(S) \leq \hat V(S)$ a.s.\\
	On the other hand, 
	\begin{eqnarray}
	\hat V(S)&=& \esssup_{\tau\in\mathcal{T}_S^p}E[V_p(\tau) {\bf 1}_{\{\alpha V_p(\tau)
		\leq \phi(\tau)\}}|\mathcal{F}_{S^-}] \\
	&\leq&   \esssup_{\tau\in\mathcal{T}_{S}^p, \tau \geq \tau^\alpha(S)}E[V_p(\tau)|\mathcal{F}_{S^-}]\\
	& \leq&\nonumber E[M_{\tau^\alpha(S)}-A_{\tau^\alpha(S)}-C_{\tau^\alpha(S)^{-}}|\mathcal{F}_{S^-}]\\ \nonumber
	& \leq &E[V_S|\mathcal{F}_{S^-}]=V(S).
	\end{eqnarray}
	Thus, these inequalities become equalities. It follows that  $$\hat V(S)= V(S)= E[M_{\tau^\alpha(S)}-A_{\tau^\alpha(S)}-C_{\tau^\alpha(S)^{-}}|\mathcal{F}_{S^-}]. $$
This implies that
	$A_{\hat \tau(S)}=A_S $ a.s. and $C_{\hat \tau(S)^-}=C_{S^-}$ a.s.	
	%\subsubsection{ Minimal optimal stopping time}
\fproof

\begin{lemma}\label{lem}	Suppose the reward $(\phi(\tau), \tau\in \mathcal{T}_0^p)$ is right limited (RL) and left limited (LL), right USCE  and $V_P(0)< \infty$. Let $S$ $\in$ $\mathcal{T}_0^p$. We define  the limiting stopping time  $\hat \tau(S)$ by 
	\begin{equation*}\label{hattheta}
	\hat \tau(S):= \lim_{\alpha\uparrow 1} \uparrow \tau^\alpha(S)\;\;\mbox{a.s}.
	\end{equation*} 
	Let $$H_S^-=\{\tau^\alpha(S)\;\; \mbox{increases strictly to}\;\; \hat\tau(S) \},$$
	and 
	$$H_S= (H_S^-)^c\cap \{ V_p(\hat \tau(S))= \phi(\hat\tau(S)) \},$$
	and 
	$$H_S^+= (H_S^-)^c \cap \{ V_p(\hat \tau(S))> \phi(\hat\tau(S)) \}.$$
	Then, for each predictable time $S \in \mathcal{T}_0^p$, 
	\begin{eqnarray}
	V_p(S)&=&  E [\phi(\hat \tau(S)^-)1_{H_S^-}+\phi(\hat \tau(S))1_{H_S}+\phi(\hat \tau(S)^+)1_{H_S^+}|\mathcal{F}_{S^-}]. 
	\end{eqnarray} 
	This is equivalent to 
	$$\sup_{\tau \in \stops}E[\phi(\tau)]= E [\phi(\hat \tau(S)^-)1_{H_S^-}+\phi(\hat \tau(S))1_{H_S}+\phi(\hat \tau(S)^+)1_{H_S^+}].$$\end{lemma}
\dproof  
	To simplify the notation, in the following, the stopping time $\tau^{\alpha}(S)$ will be denoted by $\tau^{\alpha}$ and $\hat \tau(S)$ by $\hat \tau$ .
	An application of  Lemma \ref{stepun} yields 
	\begin{equation*}
	\alpha V_p(\tau^{\alpha}(S) )\leq  \phi(\tau^{\alpha}(S) ) \,\,\,{\rm a.s.}
	\end{equation*}
	for each 
	$\alpha \in ]0, 1[.$
	If we denote also $V$ the process which aggregates the family $V_p$, \\
	On $H^-$,
	$$\phi(\hat \tau^-)=V(\hat \tau^-)\;\mbox{a.s  and } \quad \Delta A_{\hat \tau}=0.$$
	On $H^+$, 
	$$\phi(\hat \tau^+)=V(\hat \tau^+)\;\mbox{a.s  and } \quad \Delta C_{\hat \tau}=0. $$
	We have by Lemma \ref{lemeq} ,
	on  $H^-\;$, 
	$A_{\hat \tau}=A_S=A_{\hat \tau^-}\quad $, $C_{\hat \tau^-}=C_{S^-}$.\\
	On  $H$ and $H^+$, $C_{\hat \tau}=C_{\hat \tau^-}=C_{S^-}$.\\
	On $H^+$,	$C_{\hat \tau}=C_{\hat \tau^-}$, thus
	\begin{eqnarray*}
		V(S)= V_S &=& E [(M_{\hat \tau^-}-A_{\hat \tau^-}-C_{\hat \tau^-})1_{H^-}|\mathcal{F}_{S^-}]+ E [(M_{\hat \tau^-}-A_{\hat \tau}-C_{\hat \tau^-})1_{H}|\mathcal{F}_{S^-}]\\
		&+& E [(M_{\hat \tau}-A_{\hat \tau}-C_{\hat \tau})1_{H^+}|\mathcal{F}_{S^-}]\\ \nonumber
		&=&  E [V_{\hat \tau^-}1_{H^-}+V_{\hat \tau}1_{H}+V_{\hat \tau^+}1_{H^+}|\mathcal{F}_{S^-}] \\ \nonumber 
		&=&  E [V(\hat \tau^-)1_{H^-}+V(\hat \tau)1_{H}+V(\hat \tau^+)1_{H^+}|\mathcal{F}_{S^-}] \\ \nonumber 
		&=&  E [\phi(\hat \tau^-)1_{H^-}+\phi(\hat \tau)1_{H}+\phi(\hat \tau^+)1_{H^+}|\mathcal{F}_{S^-}]. 
	\end{eqnarray*} 
Hence the result.
\fproof

\end{spacing}

\end{document}